\documentclass[11pt,leqno,oneside]{article}

\usepackage{amsmath}
\usepackage{amsthm}
\usepackage{amstext}
\usepackage{amsopn}
\usepackage{amsbsy}

\newtheorem*{main-theorem}{Main Theorem}
\newtheorem{theorem}{Theorem}
\newtheorem{corollary}{Corollary}
\newtheorem{lemma}{Lemma}

\begin{document}
\pagestyle{plain}

\title{Estimates on the Lower Bound of the First Gap
\thanks{2000 Mathematics Subject Classification Primary 35J10; Secondary 35P15, 53C21}}

\author{Jun LING}

\date{}
\maketitle

\begin{abstract}
We give a new lower bound for the first gap $\lambda _2 -
\lambda_1$ of the Dirichlet eigenvalues  of the Schr{\"o}dinger
operator on a bounded convex domain $\Omega$ in R$^n$ or S$^n$ and
greatly sharpens the previous estimates. The new bound is explicit
and computable.
\end{abstract}

%%%%%%%%%%%%%%%%%%Introduction%%%%%%%%%%%%%%%%%%
\section{Introduction}\label{sec-intro}
In this paper, we give a new estimate on the lower bound of the
gap of the first two Dirichlet eigenvalues  of the Schr{\"o}dinger
operator on a bounded strictly convex domain $\Omega $ in R$^n$ or
S$^n$. Let $\lambda_1$ and $\lambda_2$ be the first two Dirichlet
eigenvalues of the Schr{\"o}dinger operator $-\Delta+V$ to the
eigenvalue problem
\begin{equation}                                         \label{basic1}
-\Delta u + Vu =\lambda u \quad \textrm{in }\Omega, \qquad
u=0\quad\textrm{on } \partial\,\Omega,
\end{equation}
where $\Delta$ is the Laplacian on $\Omega$, $V:
\Omega\rightarrow\textrm{R}^1$ a nonnegative convex smooth
function, $\Omega$ is a bounded domain in R$^n$ or S$^n$ whose
second fundamental form of the boundary with respect to the
outward normal to the boundary is positive definite. It is an
interesting and important problem to find a lower bound for the
first gap $\lambda_2-\lambda_1>0$. There has been a lot of work on
this problem. See \cite{a} and \cite{Ban-MenHer} for the
references on the problem. In 1983, M.~van~den~Berg \cite{be}
conjectured that the lower bound is $3\pi^2/d^2$, where $d$ is the
diameter of the domain. See also S.~T.~Yau's \textit{Problem
Section} in \cite{sy} and M.~Ashbaugh \cite{a} about the
conjecture. In recent years, M.~Ashbaugh and R.~Benguria
\cite{ab1}\cite{ab2} and R. Ba\~{n}uelos and P.J.
M\'{e}ndez-Hern\'{a}ndez \cite{Ban-MenHer}, R.~Ba\~{n}uelos and
P.~Kroger \cite{Ban-Kro}, and B.~Davis \cite{Davis} proved the
conjecture for some special potential $V$ and for some special
class of symmetric domains in R$^2$. For a general bounded convex
domain in R$^n$, Singer, Wong, Yau and Yau \cite{swyy} showed that
$\lambda_2-\lambda_1\geq \pi^2/(4d^2)$.  Yu and Zhong \cite{yz2}
later removed the factor $4$, using the interior log-convexity of
a positive first eigenfunction. Lee and Wang \cite{lw} showed that
one still has interior log-convexity if the domain is in S$^n$ and
therefore the above estimate remains true for the Laplacian on a
bounded convex domain in S$^n$. The author \cite{ling1} proved
that global log-convexity holds if the domain in R$^n$ or S$^n$ is
strictly convex and therefore one has the strict lower bound $
\lambda_2-\lambda_1>\pi^2/d^2$. R.~G.~Smits \cite{Smits} gave an
alternative derivation of the last inequality. In this paper, we
give a new general bound for a general bounded convex domain. The
new bound is explicit and computable.

In \cite{yau4}, S.~T.~Yau gave a estimate $ \lambda_2-\lambda_1
\geq \theta\,\pi^2/d^2+ 2\cos^2(\sqrt{\theta}\pi)\,\alpha $, where
$\theta$ is any constant with $0\leq \theta\leq 1/4$, and
$\alpha>0$ is a quantity determined by the global log-convexity
and is defined in (\ref{alpha-def}). Let us first recall the
notion of "log-convexity". Let $f$ be a positive first
eigenfunction. Let $\alpha(x)=\inf_ {\tau\in\textrm{T}_x\Omega,
|\tau|=1} \left[\nabla^2\left(-\ln f\right)\right] \left(\tau,
\tau\right)\left(x\right).$ It is known from the work of Brascamp
and Lieb \cite{bl}, Caffarelli and Friedman \cite{cf}, Korevaar
\cite{k}, Korevaar and Lewis \cite{kl}, and Lee and Wong \cite{lw}
that the interior log-convexity $\alpha (x)
> 0$ for $x$ in $\Omega$ holds. Let $\alpha$ be the constant
\begin{equation}                                \label{alpha-def}
\alpha=\inf_{x\in\Omega} \alpha(x).
\end{equation}
The author \cite{ling1} showed that the global log-convexity
$\alpha
> 0$ holds.
Yau \cite{yau4} gave an interesting estimate on the lower bound
$\alpha$ in terms of the potential $V$.

Our main result is the following Theorem \ref{main-thm}.

\begin{theorem}                                    \label{main-thm}
If $\Delta$ is the Laplacian in \textup{R}$^n$ or \textup{S}$^n$
and if $\lambda_1$ and $\lambda_2$ are the first two Dirichlet
eigenvalues of the Schr{\"o}dinger operator $-\Delta + V$ with a
nonnegative convex potential $V$ on a bounded strictly convex
domain $\Omega$ in \textup{R}$^n$ or \textup{S}$^n$, then we have
the following estimate
\[
\lambda_2 -\lambda_1 \geq \frac{\pi^2}{d^2} + \frac{31}{50}\alpha,
\]
where $\alpha$ is the quantity in (\ref{alpha-def}) and $d$ is the
diameter of $\Omega$.
\end{theorem}

\noindent \textit{Remark 1.} If the domain $\Omega$ has a certain
symmetric property that the "midrange" of the ratio of
eigenfunctions is zero, then Theorem \ref{thm4} shows
\[
\lambda_2 -\lambda_1 \geq \frac{\pi^2}{d^2} + \alpha.
\]

 \noindent \textit{Remark 2.} The above result holds when
$\Omega$ is in a manifold with non-negative Ricci curvature and
positive $\alpha$.

In order to improve the known results on the gap estimate, we need
to construct suitable test functions where detailed technical work
is essential. In the last section we construct the test function
$\xi$. We explore the properties of the function $\xi$, the
Zhong-Yang function $\eta$ and the ratio $\xi/\eta$. Those
properties are essential to the construction of the suitable test
functions. Because those functions are complicated combinations of
trigonometric and rational functions, the needed properties such
as monotonic and convex properties are hard to prove. In the past,
though we know that many nice properties might be true, only a few
of them could be proven strictly in mathematics by the canonical
calculus method. We are able to prove those properties effectively
now by studying the differential equations those functions
satisfied and using the Maximum Principle. Since the constructions
and proofs in that part are quite technical by nature, we put them
in the last section. Readers may refer to that section when in
need. The functions $\xi$ and $\eta$ and their properties have
also other important applications. In our recent estimates on the
first non-zero eigenvalues of compact manifolds, function $\xi$
plays major role.  We prove the main result in Section
\ref{sec-proofs}. Last section is for deriving some preliminary
estimates and the conditions for test functions and for proving
the properties of the functions used in the proof of the Theorem
\ref{main-thm}.

%%%%%%%%%Proof of Main Teorems%%%%%%%%%%%%%%%%%%%%
\section{Proof of the Main Result}\label{sec-proofs}
Let $f_2$ be the second eigenfunction and $f$ a positive first
eigenfunction of Problem (\ref{basic1}). It is known $f_2/f$
changes its signs in $\Omega$ (see \cite{cha}) and is smooth up to
the boundary $\partial\Omega$ (see \cite{swyy}) and can be
normalized so that
\[
\sup_{\Omega}f_2/f=1,
\quad \inf_{\Omega}f_2/f=-k, \quad \textrm{and}\quad 0<k\leq 1.
\]
Let $\lambda=\lambda_2-\lambda_1$,
\begin{equation}                                \label{v-def}
v=[ f_2/f - (1-k)/2 ]/ [ (1+k)/2 ].
\end{equation}
 Then the function $v$ satisfies the following
\begin{equation}
\Delta v=-\lambda (v+a)-2\nabla v\,\nabla (\ln f)  \quad
\textrm{in }\Omega,                             \label{basic2}
\end{equation}
\begin{equation}
\frac{\partial v}{\partial N}=0  \quad \textrm{on }\partial
\Omega, \label{basic2.001}
\end{equation}
\begin{equation} \label{basic2.01}
\sup_{\Omega}v=1\quad \textrm{and} \quad \inf_{\Omega}v=-1
\end{equation}
where $N$ is the the outward normal of $\partial \Omega$, and
\begin{equation}
a=(1-k)/(1+k)                  \label{a-def}
\end{equation}
is the "midrange" of the ratio $f_2/f$. Note that $0\leq a<1$.

We set in this paper, unless otherwise stated,
\begin{equation}                        \label{delta-def}
\lambda=\lambda_2-\lambda_1\quad \textrm{and}\quad
\delta=\alpha/\lambda.
\end{equation}
and let
\[ Z(t)=\max_{x\in\bar{\Omega},t=\sin^{-1} \left(v\left(x\right)/b\right)}
\frac{\left |\nabla v\right |^2}{b^2-v^2}/\lambda.
\]
for $t\in [-\sin^{-1}(1/b), \sin^{-1}(1/b)]$.

We have the following  estimates (\ref{delta-bound}) and
(\ref{basic5}), Lemma \ref{barrier-lemma} and corollaries
\ref{corollary1} and \ref{corollary2}. The proofs are in the last
section.

\begin{equation}                \label{delta-bound}
0<\delta \leq \frac12.
\end{equation}

\begin{equation} \label{basic5}
Z(t)\leq 1+ a\quad t\in [-\sin^{-1}(1/b),  \sin^{-1}(1/b)].
\end{equation}

\begin{lemma}                                               \label{barrier-lemma}
If the function $z:\ [-\sin^{-1} (1/b),\,\sin^{-1} (1/b)]\mapsto
\mathbf{R}^1$ satisfies the following
\begin{enumerate}
 \item $z(t)\geq Z(t) \qquad t\in [-\sin^{-1}(1/b),  \sin^{-1}(1/b)]$,
 \item there exists some $x_0\in \bar{\Omega}$
       such that at point $t_0=\sin^{-1} (v(x_0)/b)$ \linebreak
       $z(t_0)=Z(t_0)$,
 \item $z(t_0)>0$,
\end{enumerate}
then we have the following
\begin{eqnarray}
0&\leq&\frac12z''(t_0)\cos^2t_0 -z'(t_0)\cos t_0\sin t_0 - z(t_0)
            +1+c\sin t_0 -2\delta \cos^2t_0\nonumber\\
 & & {}-\frac{z'(t_0)}{4z(t_0)}\cos t_0[z'(t_0)\cos t_0 -2z(t_0)\sin t_0
            + 2\sin t_0 + 2c].\label{barrier-eq}
\end{eqnarray}
\end{lemma}

\begin{corollary}                                                   \label{corollary1}
If in addition to the above conditions 1-3 in Lemma
\ref{barrier-lemma}, $z'(t_0)\geq 0$ and $1-c\leq z(t_0)\leq 1+a$,
then we have the following
\begin{equation}
0\leq\frac12z''(t_0)\cos^2t_0 -z'(t_0)\cos t_0\sin t_0-z(t_0) +
1+c\sin t_0-2\delta \cos^2t_0.\nonumber
\end{equation}
\end{corollary}

\begin{corollary}                                          \label{corollary2}
If $a=0$, which is defined in (\ref{a-def}), and if  in addition
to the above conditions 1-3 in Lemma \ref{barrier-lemma},
$z'(t_0)\sin t_0\geq 0$ and $z(t_0)\leq 1$, then we have the
following
\begin{equation}
0\leq\frac12z''(t_0)\cos^2t_0 -z'(t_0)\cos t_0\sin t_0 -z(t_0)+
1-2\delta \cos^2t_0.\nonumber
\end{equation}
\end{corollary}

We now prove our result.

\begin{theorem}       \label{thm3}
If $a>0$ and $\mu\delta \leq \frac{4}{\pi^2}a$ for a constant
$\mu\in (0,1]$, then
\[
\lambda_2 -\lambda_1 \geq \frac{\pi^2}{d^2} + \mu\alpha
\]
\end{theorem}

\begin{proof}\quad Let $\mu_{\epsilon}=\mu-\epsilon>0$ for small
positive constant $\epsilon$. Take $b>1$ close to $1$ such that
$\mu_{\epsilon}\delta < \frac{4}{\pi^2}c$. Let
\begin{equation}                                \label{z-def}
z(t)=1+c\eta(t) +\mu_{\epsilon}\delta\xi(t),
\end{equation}
where $\xi$ and $\eta$ are the functions defined by (\ref{xi-def})
and (\ref{eta-def}), respectively.  Let
$\bar{I}=[-\sin^{-1}(1/b),\sin^{-1}(1/b)]$. We claim that
\begin{equation}\label{4.1}
Z(t)\leq z(t)\qquad \textrm{for }t\in \bar{I}.
\end{equation}
By Lemma \ref{xi-lemma} and Lemma \ref{eta-lemma} we have
\begin{eqnarray}
& &{}\frac{1}{2}z''\cos ^2t-z'\cos t\sin t-z
    =-1-c\sin t+ 2\mu_{\epsilon}\delta\cos^2t,          \label{z-eq}\\
& &{}z'(t)> 0\label{z'-geq0}\\
& &{}0<1-\frac{a}{b}=z(-\frac{\pi}{2})\leq z(t)\leq
z(\frac{\pi}{2})=1+\frac{a}{b}\leq 1+a,          \label{z-endpts}
\end{eqnarray}
where (\ref{z'-geq0}) is due to the following.
\begin{eqnarray}
z'(t)=c\eta'(t)+\mu_{\epsilon}\delta\xi'(t)&=&\mu_{\epsilon}\delta\eta'(t)\left(
\frac{c}{\mu_{\epsilon}\delta}+\frac{\xi'(t)}{\eta'(t)}\right)\nonumber\\
&\geq & \mu_{\epsilon}\delta
\eta'(t)(\frac{c}{\mu_{\epsilon}\delta}-\frac{\pi^2}{4}) >
0.\nonumber
\end{eqnarray}
Let $P\in\mathbf{R}^1$ and $t_0\in
[-\sin^{-1}(1/b),\sin^{-1}(1/b)]$ such that
\[ P=\max_{t\in \bar{I}}\left(Z(t)-z(t)\right)=Z(t_0)-z(t_0).
\]
Thus
\begin{equation}\label{4.2}
Z(t)\leq z(t)+P\quad \textrm{for }t\in
\bar{I}\qquad\textrm{and}\qquad Z(t_0)=z(t_0)+P.
\end{equation}
Suppose that $P>0$ Then $z+P$ satisfies the inequality in
Corollary \ref{corollary1} of Lemma \ref{barrier-lemma}. Then
\begin{eqnarray}
&{}&z(t_0)+P=Z(t_0)\nonumber\\
&{}&\leq  \frac12(z+P)''(t_0)\cos^2 t_0-(z+P)'(t_0)\cos t_0
 \sin t_0+1+c\sin t_0-2\delta \cos^2 t_0\nonumber\\
&{}&=\frac12z''(t_0)\cos^2t_0-z'(t_0)\cos t_0\sin
t_0+1+c\sin t_0-2\delta \cos^2 t_0\nonumber\\
&{}&\leq\frac12z''(t_0)\cos^2t_0-z'(t_0)\cos t_0\sin
t_0+1+c\sin t_0-2\mu_{\epsilon}\delta \cos^2 t_0\nonumber\\
&{}&=z(t_0).\nonumber
\end{eqnarray}
This contradicts the assumption $P>0$. Thus $P\leq 0$ and
(\ref{4.1}) must hold. Now we have
\[
|\nabla t|^2\leq\lambda z(t) \qquad \textrm{for } t\in\bar{I},
\]
that is
\begin{equation}\label{4.3}
\sqrt{\lambda}\geq \frac{|\nabla t|}{\sqrt{z(t)}}.
\end{equation}
Let $q_1$ and $q_2$ be two points in $\bar{\Omega}$ such that
$v(q_1)=-1$ and $v(q_2)=1$ and let $L$ be the minimum geodesic
segment between $q_1$ and $q_2$. $L$ lies on $\bar{\Omega}$
completely, since $\bar{\Omega}$ is convex. We integrate both
sides of (\ref{4.3}) along $L$ and change variable and let
$b\rightarrow 1$. Then
\begin{equation}\label{4.4}
\sqrt{\lambda}d\geq \int_{L}\,\frac{|\nabla
t|}{\sqrt{z(t)}}dl=\int_{-\frac{\pi}{2}}^{\frac{\pi}{2}}
\frac{1}{\sqrt{z(t)}}\,dt \geq \frac{\left(\int_{-\pi/2}^{\pi/2}\
\,dt\right)^\frac32}{(\int_{-\pi/2}^{\pi/2}\ z(t)\,dt)^{\frac12}}
\geq \left( \frac{\pi^3}{\int_{-\pi/2}^{\pi/2}\  z(t)\,dt}
\right)^{\frac12}.
\end{equation}
Square the two sides. Then
\[
\lambda \geq \frac{\pi^3}{d^2\int_{-\pi/2}^{\pi/2} \ z(t)\,dt}.
\]
Now
\[
\int_{-\frac{\pi}{2}}^{\frac{\pi}{2}}\
z(t)\,dt=\int_{-\frac{\pi}{2}}^{\frac{\pi}{2}}\ [1+ a\eta(t)+
\mu_{\epsilon}\delta \xi(t)]\,dt=(1-\mu_{\epsilon}\delta)\pi,
\]
where we used the facts that
$\int_{-\frac{\pi}{2}}^{\frac{\pi}{2}}\ \eta(t)\,dt=0$ since
$\eta$ is an even function, and that $
\int_{-\frac{\pi}{2}}^{\frac{\pi}{2}}\ \xi(t)\,dt=-\pi$ (by
(\ref{xi-int}) in the Lemma \ref{xi-lemma}). Therefore
\[
\lambda \geq
\frac{\pi^2}{(1-\mu_{\epsilon}\delta)d^2}\quad\textrm{and}\quad
\lambda \geq \frac{\pi^2}{d^2} + \mu_{\epsilon}\alpha.
\]
Letting $\epsilon\rightarrow 0$, we get
\[
\lambda \geq \frac{\pi^2}{(1-\mu\delta)d^2}\quad\textrm{and}\quad
\lambda \geq \frac{\pi^2}{d^2} + \mu\alpha.
\]
\end{proof}

\begin{theorem}                                          \label{thm4}
If the "midrange" $a=0$, then
\begin{equation}                        \label{ling1-4}
\lambda_2 -\lambda_1 \geq \frac{\pi^2}{d^2}+\alpha.
\end{equation}
\end{theorem}

\begin{proof}\quad Let
\[
y(t)=1+\delta\xi.
\]
By Lemma \ref{xi-lemma},
 for $-\frac{\pi}{2}< t<\frac{\pi}{2}$, we have
\begin{eqnarray}
& &{}\frac{1}{2}y''\cos ^2t-y'\cos t\sin t-y
    =-1+2\delta\cos^2t,\label{21.1}\\
& &{}y'(t)\sin t\geq 0,\qquad \textrm{and} \label{y'-geq0}\\
& &{}y(\pm \frac{\pi}{2})=1 \  \textrm{and }\,
0<y(t)<1.\label{23.3}
\end{eqnarray}
We need only show that $Z(t)\leq y(t)$ on $[-\pi/2,\pi/2]$. If it
is not true, then there is $t_0$ and a number $P>0$ such that
$P=Z(t_0)-y(t_0)=\max Z(t)-y(t)$. Note that $y(t)+P\geq
1-\frac12(\frac{\pi^2}{4}-1)+P>0$. So $y+P$ satisfies the
inequality in the Corollary \ref{corollary2} in Lemma
\ref{barrier-lemma}. Therefore
\begin{eqnarray}
&{}&y(t_0)+P=Z(t_0)\nonumber\\
&{}&\leq  \frac12(y+P)''(t_0)\cos^2 t_0-(y+P)'(t_0)\cos t_0
 \sin t_0+1-2\delta \cos^2 t_0\nonumber\\
&{}&=\frac12y''(t_0)\cos^2t_0-y'(t_0)\cos t_0\sin
t_0+1-2\delta \cos^2 t_0\nonumber\\
&{}&=y(t_0).\nonumber
\end{eqnarray}
This contradicts the assumption $P>0$. The rest of the proof is
similar to that of Theorem \ref{thm3}, just noticing that
$\delta\leq \frac12<\frac{4}{\pi^2-4}$.
\end{proof}

\begin{proof}[Proof of Theorem \ref{main-thm}]\quad
Since $0\leq a<1$, either $a=0$ or $0<a<1$.

If $a=0$, then we apply Theorem \ref{thm4} to get the bound with
$\mu=1$,
\[
\lambda \geq \frac{\pi^2}{d^2} + \alpha.
\]

If $0<a<1$, then there are several cases altogether.

\begin{itemize}
\item(I): $\ $$a\geq\frac{\pi^2}{4}\delta$. \item(II):
$a<\frac{\pi^2}{4}\delta$.
    \begin{itemize}
    \item(II-a): $a\geq 0.765$.
    \item(II-b): $0<a<0.765$.
        \begin{itemize}
        \item(II-b-1): $a\geq 1.53\delta$.
        \item(II-b-2): $a<1.53\delta$.
        \end{itemize}
    \end{itemize}
\end{itemize}

For Case (I): $\ $ $0<a<1$ and $a\geq\frac{\pi^2}{4}\delta$, we
apply Theorem \ref{thm4} for $\mu=1$ to get the following lower
bound
\[
\frac{\pi^2}{d^2} + \alpha
\]

For Case (II-a): $0.765\leq a<\frac{\pi^2}{4}\delta$, we apply
Theorem \ref{thm3} with $\mu=\frac{4}{\pi^2}\frac{a}{\delta}$
since $(\frac{4}{\pi^2}\frac{a}{\delta})\,\delta\leq
\frac{4}{\pi^2}a$ and $0<\frac{4}{\pi^2}\frac{a}{\delta}<1$. Then
\[
\lambda\geq
\frac{\pi^2}{d^2}+\frac{4}{\pi^2}\frac{a}{\delta}\alpha
=\frac{\pi^2}{d^2}+\frac{4a}{\pi^2}\lambda
\]
Thus
\[
\lambda \geq \frac{1}{1-\frac{4a}{\pi^2}}\,\frac{\pi^2}{d^2}.
\]
On the other hand we have bound (\ref{delta-bound}),
\[
\lambda \geq 2\alpha.
\]
The above two estimates give
\[
\lambda \geq \frac{\pi^2}{d^2} +\frac{4a}{\pi^2}\,2\alpha \geq
\frac{\pi^2}{d^2} +\frac{8(0.765)}{\pi^2}\alpha
>
\frac{\pi^2}{d^2}+ \frac{31}{50}\alpha.
\]
The theorem is proved in this case.

For Case (II-b-1): $0<a<0.765$, $a<\frac{\pi^2}{4}\delta$ and
$a\geq 1.53\delta$, we apply Theorem \ref{thm3} with with
$\mu=\frac{4}{\pi^2}\frac{a}{\delta}$ since
$(\frac{4}{\pi^2}\frac{a}{\delta})\,\delta\leq \frac{4}{\pi^2}a$
and $0<\frac{4}{\pi^2}\frac{a}{\delta}<1$. Then
\[
\lambda\geq
\frac{\pi^2}{d^2}+\frac{4}{\pi^2}\frac{a}{\delta}\alpha \geq
\frac{\pi^2}{d^2}+\frac{4}{\pi^2}\frac{153}{100}\alpha
>\frac{\pi^2}{d^2}+\frac{31}{50}\alpha,
\]
which is what we want to prove.

For the remaining Case (II-b-2): $0<a<0.765$,
$a<\frac{\pi^2}{4}\delta$ and $a<1.53\delta$, we define a function
$z$ by
\[
z(t)=1+c\eta(t)+(\delta-\sigma c^2)\xi(t) \qquad\textrm{on
}[-\sin^{-1}\frac1b, \sin^{-1}\frac1b],
\]
where
\begin{equation}                                      \label{sigma-def}
\sigma =\frac{\tau}{\left([\,\frac32-\frac{\pi^2}{8}-
(\frac{\pi^2}{32}-\frac16)\frac{153}{100}]\frac{200}{153}-
\frac{(\frac{8}{3\pi}-\frac{\pi}{4})^2}{[-1
+(12-\pi^2)\frac{100}{153}]}\right)c}
\end{equation}
and
\begin{equation}                                      \label{tau-def}
\tau =
\frac{2}{3\pi^2}\left(\frac{4}{3(4-\pi)}+\frac{3(4-\pi)}{4}-2\right).
\end{equation}

Let $\bar{I}=[-\sin^{-1}\frac1b, \sin^{-1}\frac1b]$. We now show
that
\begin{equation}\label{Z-leq-z}
Z(t)\leq z(t)\qquad\textrm{on }\bar{I}.
\end{equation}

If (\ref{Z-leq-z}) is not true, then there exists a constant $P>0$
and $t_0$ such that
\[
Pc^2=\frac{Z(t_0)-z(t_0)}{-\xi(t_0)}=\max_{t\in[-\sin^{-1}\frac1b,
\sin^{-1}\frac1b]}\frac{Z(t)-z(t)}{-\xi(t)}.
\]
Let $w(t)=z(t)-Pc^2\xi(t)=1+c\eta(t) +m\xi(t)$, where $m =
\delta-\sigma c^2 -Pc^2$. Then
\[
Z(t)\leq w(t)\quad\textrm{on }\bar{I}\qquad \textrm{and}\qquad
Z(t_0)=w(t_0).
\]
By Lemma \ref{z-positive-lemma}, $w(t_0)>0$. So $w$ satisfies
(\ref{barrier-eq}) in Lemma \ref{barrier-lemma},
\[
0\leq -2(\sigma + P)c^2 \cos^2t_0
  -\frac{w'(t_0)}{4w(t_0)}\cos t_0\left(\frac{8c}{\pi}\cos t +4m t\cos
  t\right).
\]
We used (\ref{xi-eq}), (\ref{xi-eq2}), (\ref{eta-eq}) and
(\ref{eta-eq2}) to get the above inequality. Thus
\begin{equation}\label{48.0}
\sigma + P\leq
  -\frac{w'(t_0)}{2c^2w(t_0)}\left(\frac{2c}{\pi} +m t\right)
  =-\frac{\eta'(t_0)}{\pi
  w(t_0)}\left(1+\frac{m\xi'(t_0)}{c\eta'(t_0)}\right)
     \left(1 +\frac{\pi m}{2c}t_0\right).
\end{equation}
The righthand side is not positive for $t_0\geq 0$, by Lemmas
\ref{xi-lemma} and \ref{eta-lemma}. Thus $t_0<0$, and
\begin{eqnarray}
 {}&{}&-\left(1+\frac{m\xi'(t_0)}{c\eta'(t_0)}\right)
     \left(1 +\frac{\pi m}{2c}t_0\right)
\nonumber\\
{}&{}&=\frac{2\xi'(t_0)}{\pi t_0\eta'(t_0)}\left(\frac{\pi
t_0\eta'(t_0)}{2\xi'(t_0)}+\frac{\pi m}{2c}t_0\right)\left(-1
-\frac{\pi m}{2c}t_0\right)\nonumber\\
{}&{}&\leq \frac14\frac{2\xi'(t_0)}{\pi t_0\eta'(t_0)}\left(
\frac{\pi
t_0\eta'(t_0)}{2\xi'(t_0)}-1 \right)^2\nonumber\\
{}&{}&= \frac14\left( \frac{2\xi'(t_0)}{\pi
t_0\eta'(t_0)}+(\frac{2\xi'(t_0)}{\pi
t_0\eta'(t_0)})^{-1}-2\right). \nonumber
\end{eqnarray}
By Lemmas \ref{xi-lemma} and \ref{eta-lemma}, we have
$2(3-\frac{\pi^2}{4})\leq \frac{\xi'(t)}{t}\leq \frac43$ and
$2(\frac{4}{\pi}-1)\leq \eta'(t)\leq \frac{8}{3\pi}$. So
\[
\frac{3(12-\pi^2)}{8}\leq\frac{2\xi'(t_0)}{\pi t_0\eta'(t_0)}\leq
\frac{4}{3(4-\pi)}.
\]
Note that the function $f(t)=t+\frac1t-2$ achieves it maximum on
$[A, B]$ not containing $0$ at an endpoint. Therefore
\[
\left|-\left(1+\frac{m\xi'(t_0)}{c\eta'(t_0)}\right)
     \left(1 +\frac{\pi m}{2c}t_0\right)\right|\leq
     \frac14\left(\frac{4}{3(4-\pi)}+\frac{3}{3(4-\pi)}-2\right).
\]

Now (\ref{48.0}) becomes
\begin{equation}  \label{49}
\sigma + P\leq \frac{\tau}{w(t_0)}.
\end{equation}

On the other hand, by Lemma \ref{z-positive-lemma},
\begin{equation}\label{50}
z(t_0)\geq
\left([\,\frac32-\frac{\pi^2}{8}-(\frac{\pi^2}{32}-\frac16)\frac{153}{100}]\frac{200}{153}-
\frac{(\frac{8}{3\pi}-\frac{\pi}{4})^2}{[-1
+(12-\pi^2)\frac{100}{153}]}\right)c=\frac{\tau}{\sigma}>0.
\end{equation}
Since $-P\xi(t_0)\geq 0$, we have $w(t_0)\geq z(t_0)$. This fact,
(\ref{49}) and (\ref{50}) imply that for $P>0$
\[
\sigma + P<\sigma,
\]
which is impossible.

Therefore we have the estimate (\ref{Z-leq-z}). Now we proceed as
in the proof of Theorem \ref{thm3}. We get the following
\[
\lambda d^2\geq \frac{\pi^3}{\pi[1-(\delta-\sigma c^2)]}.
\]
Since $\delta -\sigma c^2>0.625\delta$ by Lemma
\ref{z-positive-lemma}, we have
\[
\lambda \geq \frac{1}{[1-(\delta-\sigma
c^2)]}\frac{\pi^2}{d^2}\geq
\frac{1}{[1-0.625\delta]}\frac{\pi^2}{d^2}
\]
and
\[
\lambda \geq
\frac{\pi^2}{d^2}+0.625\alpha>\frac{\pi^2}{d^2}+\frac{31}{50}.
\]
\end{proof}

We now present a Lemma that is used in the proof of the Theorem
\ref{main-thm}.

\begin{lemma}\label{z-positive-lemma}
If $a<1.53\delta$ and $0<a<0.765$ then
\[
z(t)=1+c\eta(t)+\delta\xi(t)
\]
\[
\geq
\left([\,\frac32-\frac{\pi^2}{8}-(\frac{\pi^2}{32}-\frac16)\frac{153}{100}]\frac{200}{153}-
\frac{(\frac{8}{3\pi}-\frac{\pi}{4})^2}{[-1
+(12-\pi^2)\frac{100}{153}]}\right)c>0,
\]
for $t\in[-\pi/2, \pi/2]$ and
\[
\delta-\sigma c^2 \approx 0.625162283437>0.625\delta,
\]
where $c=a/b$ and $b>1$ is any constant and $\sigma$ is the
constant in (\ref{sigma-def}).
\end{lemma}
\begin{proof}\quad By Lemmas \ref{r-lemma}, Lemma \ref{xi-lemma} and \ref{eta-lemma},
the function $z$ on $[-\pi/2, \pi/2]$ has a unique critical point
$t_1\in (-\pi/2, 0)$  if $0<a<\frac{\pi^2}{4}\delta$ and $z$ is
decreasing on $[-\pi/2, t_1]$ and increasing on $[t_1, \pi/2]$.
Therefore
\[
\min_{[-\pi/2, \pi/2]}z=\min_{[-\pi/2, 0]}z=z(t_1).
\]
So we need only consider the restricted function $z|_{[-\pi/2,
0]}$ for the minimum.

Now first consider the Taylor expansion of $\xi$ at $0$ for
$t\in[-\pi/2, 0]$. By Lemma \ref{xi-lemma},
$\xi(0)=-\frac{\pi^2}{4}+1$, $\xi'(0)=0$ and
$\xi''(0)=2(3-\frac{\pi^2}{4})$ and $\xi'''(t)<0$ on $(-\pi/2,
0)$.

Thus

\[
\xi(t)=\xi(0)
+\xi'(0)+\frac{\xi''(0)}{2!}t^2+\frac{\xi'''(t_2)}{2!}t^3
\]
\[
\geq\xi(0) +\xi'(0)+\frac{\xi''(0)}{2!}t^2
\]
\[
=-(\frac{\pi^2}{4}-1)+(3-\frac{\pi^2}{4})t^2,
\]
where $t_2$ is a constant in $(t, 0)$. Similarly, using the data
$\eta(-\pi/2)=-1$, $\eta'(-\pi/2)=\frac{8}{3\pi}$ and
$\eta'''(t)>0$ on $(-\pi/2, 0)$ (actually on $[-\pi/2, \pi/2]$),
and the Taylor expansion of $\eta$ at $-\pi/2$,  we have for $t\in
[-\pi/2, 0]$,
\[
\eta(t)=\eta(-\frac{\pi}{2})+\eta'(-\frac{\pi}{2})(t+\frac{\pi}{2})
+\frac{\eta''(-\frac{\pi}{2})}{2!}(t+\frac{\pi}{2})^2
+\frac{\eta''(t_3)}{3!}(t+\frac{\pi}{2})^3
\]
\[
\geq \eta(-\frac{\pi}{2})+\eta'(-\frac{\pi}{2})(t+\frac{\pi}{2})
+\frac{\eta''(-\frac{\pi}{2})}{2!}(t+\frac{\pi}{2})^2
\]
\[
=-1+\frac{8}{3\pi}(t+\frac{\pi}{2})-\frac{1}{4}(t+\frac{\pi}{2})^2
\]
\[=-(\frac{\pi^2}{16}-\frac13)+(\frac{8}{3\pi}-\frac{\pi}{4})t-\frac14
t^2,
\]
where $t_3$ is some constant in $(-\pi/2,t)$. Therefore on
$[-\pi/2,0]$,
\[
z(t)=1+c\eta(t)+\delta\xi(t)
\]
\[
\geq
1-(\frac{\pi^2}{16}-\frac13)c-(\frac{\pi^2}{4}-1)\delta+(\frac{8}{3\pi}-\frac{\pi}{4})ct+[-\frac14c
+(3-\frac{\pi^2}{4})\delta] t^2
\]

Let $\nu=1.53$ and $a_0=0.765$. That $a\leq \nu\delta$ implies
$c=a/b<\nu \delta$, where $b>1$ is a constant.  Using conditions
(\ref{delta-bound}) $\delta\leq \frac{n-1}{2n}<\frac12$ and $a\leq
a_0$, we get

\[
1-(\frac{\pi^2}{16}-\frac13)c-(\frac{\pi^2}{4}-1)\delta
\]
\[
\geq1-(\frac{\pi^2}{16}-\frac13)\nu\delta-(\frac{\pi^2}{4}-1)\delta
\]
\[
\geq \frac32-\frac{\pi^2}{8}-(\frac{\pi^2}{32}-\frac16)\nu
\]
\[
>
\left(\frac32-\frac{\pi^2}{8}-(\frac{\pi^2}{32}-\frac16)\nu\right)\frac{1}{a_0}c
\]

and

\[
1+c\eta(t)+\delta\xi(t)
\]
\[
\geq
\left(\frac32-\frac{\pi^2}{8}-(\frac{\pi^2}{32}-\frac16)\nu\right)\frac{1}{a_0}c
+(\frac{8}{3\pi}-\frac{\pi}{4})ct+[-\frac14c
+(3-\frac{\pi^2}{4})\frac{1}{\nu}c] t^2
\]
\[
=\left([\,\frac32-\frac{\pi^2}{8}-(\frac{\pi^2}{32}-\frac16)\nu]\frac{1}{a_0}
+(\frac{8}{3\pi}-\frac{\pi}{4})t+[-\frac14
+(3-\frac{\pi^2}{4})\frac{1}{\nu}] t^2\right)c
\]
\[
\geq
\left([\,\frac32-\frac{\pi^2}{8}-(\frac{\pi^2}{32}-\frac16)\nu]\frac{1}{a_0}-
\frac{(\frac{8}{3\pi}-\frac{\pi}{4})^2}{4[-\frac14
+(3-\frac{\pi^2}{4})\frac{1}{\nu}]}\right)c
\]
\[
\geq
\left([\,\frac32-\frac{\pi^2}{8}-(\frac{\pi^2}{32}-\frac16)\nu]\frac{1}{a_0}-
\frac{(\frac{8}{3\pi}-\frac{\pi}{4})^2}{[-1
+(12-\pi^2)\frac{1}{\nu}]}\right)c>0.5433>0.
\]
Let $\tau$ be the constant in (\ref{tau-def}). Then
\[
\sigma c^2 =\frac{\tau
c}{\left([\,\frac32-\frac{\pi^2}{8}-(\frac{\pi^2}{32}-\frac16)\nu]\frac{1}{a_0}-
\frac{(\frac{8}{3\pi}-\frac{\pi}{4})^2}{[-1
+(12-\pi^2)\frac{1}{\nu}]}\right)c},
\]
\[
\leq\frac{\tau\nu
\delta}{\left([\,\frac32-\frac{\pi^2}{8}-(\frac{\pi^2}{32}-\frac16)\nu]\frac{1}{a_0}-
\frac{(\frac{8}{3\pi}-\frac{\pi}{4})^2}{[-1
+(12-\pi^2)\frac{1}{\nu}]}\right)c}\approx 0.374837516563\delta
\]
and
\[
\delta-\sigma c^2 >0.625\delta.
\]
\end{proof}

%%%%%%%%%%%%%Estimates and Barriers%%%%%%%%%%%%%%%
\section{Some  Estimates and Lemmas}\label{sec-pre-es}

\begin{lemma}                           \label{xi-lemma}
Let
\begin{equation}                        \label{xi-def}
\xi(t)=\frac{\cos^2t+2t\sin t\cos t +t^2-\frac{\pi^2}{4}}{\cos^2t}
\qquad \textrm{on}\quad [-\frac{\pi}{2},\frac{\pi}{2}\,].
\end{equation}
Then the function $\xi$  satisfies the following
\begin{eqnarray}
& &{}\frac{1}{2}\xi''\cos ^2t-\xi'\cos t\sin t-\xi
    =2\cos^2t\quad \textrm{in }(-\frac{\pi}{2},\frac{\pi}{2}\,),          \label{xi-eq}\\
& &{}\xi'\cos t -2\xi\sin t =4t\cos t  \quad \textrm{in }(-\frac{\pi}{2},\frac{\pi}{2}\,),                     \label{xi-eq2}\\
& &{}\int_0^{\frac{\pi}{2}}\xi(t)\, dt= -\frac{\pi}{2} ,         \label{xi-int}\\
& &{}1-\frac{\pi^2}{4}=\xi(0)\leq\xi(t)\leq\xi(\pm
\frac{\pi}{2})=0\quad
\textrm{on }[-\frac{\pi}{2},\frac{\pi}{2}\,],                                \nonumber\\
& &{} \xi' \textrm{ is increasing on }
[-\frac{\pi}{2},\frac{\pi}{2}\,] \textrm{ and }
\xi'(\pm \frac{\pi}{2}) =\pm \frac{2\pi}{3},                     \nonumber\\
& &{}\xi'(t)< 0 \textrm{ on }(-\frac{\pi}{2},0)\textrm{ and \ }
\xi'(t)>0 \textrm{ on }(0,\frac{\pi}{2}\,), \nonumber \\
& &{}\xi''(\pm\frac{\pi}{2})=2, \ \xi''(0)=2(3-\frac{\pi^2}{4})
\textrm{ and \ } \xi''(t)> 0 \textrm{ on }
[-\frac{\pi}{2},\frac{\pi}{2}\,], \nonumber\\
& &{}(\frac{\xi'(t)}{t})'>0 \textrm{ on } (0,\pi/2\,)\textrm{ and
\ } 2(3-\frac{\pi^2}{4})\leq \frac{\xi'(t)}{t}\leq \frac43
\textrm{ on } [-\frac{\pi}{2},\frac{\pi}{2}\,],\nonumber \\
& &{}\xi'''(\frac{\pi}{2})=\frac{8\pi}{15}, \xi'''(t)< 0 \textrm{
on }(-\frac{\pi}{2},0) \textrm{ and \ }  \xi'''(t)>0 \textrm{ on
}(0,\frac{\pi}{2}\,). \nonumber
\end{eqnarray}
\end{lemma}
\begin{proof}\quad For convenience, let $q(t)= \xi'(t)$, i.
e.,
\begin{equation}                                       \label{q-def}
q(t) = \xi'(t) = \frac{2(2t\cos t +t^2\sin t +\cos^2 t \sin t
-\frac{\pi^2}{4}\sin t)}{\cos^3 t}.
\end{equation}
Equation (\ref{xi-eq}) and the values $\xi(\pm \frac{\pi}{2})=0$,
$\xi(0)=1-\frac{\pi^2}{4}$ and $\xi'(\pm \frac{\pi}{2}) =\pm
\frac{2\pi}{3}$ can be verified directly from (\ref{xi-def}) and
(\ref{q-def}) .  The values of $\xi''$ at $0$ and $\pm
\frac{\pi}{2}$ can be computed via (\ref{xi-eq}). By
(\ref{xi-eq2}), $(\xi(t)\cos^2 t)' =4t\cos^2 t$. Therefore
\newline $\xi(t)\cos^2 t=\int_{\frac{\pi}{2}}^t \ 4s\cos^2 s\,ds$,
and
\[
\int_{-\frac{\pi}{2}}^{\frac{\pi}{2}}\
\xi(t)\,dt=2\int_0^{\frac{\pi}{2}}\
\xi(t)\,dt=-8\int_0^{\frac{\pi}{2}}\left( \frac{1}{\cos^2(t)}
\int_t^{\frac{\pi}{2}}\ s\cos^2s\,ds\right)\,dt
\]
\[
=-8\int_0^{\frac{\pi}{2}}\left(\int_0^s\
\frac{1}{\cos^2(t)}\,dt\right)\ s\cos^2s\,ds
=-8\int_0^{\frac{\pi}{2}}\ s\cos s\sin s\,ds=-\pi.
\]
It is easy to see that $q$ and $q'$ satisfy the following
equations
\begin{equation}                                         \label{q-eq}
\frac12 q''\cos t -2q'\sin t -2q\cos t = -4 \sin t,
\end{equation}
and
\begin{equation}                                       \label{q'-eq}
\frac{\cos^2 t}{2(1+\cos^2 t)}(q')''-\frac{2\cos t\sin t}{1+\cos^2
t}(q')'-2(q')=-\frac{4}{1+\cos^2 t}.
\end{equation}
The last equation implies $q'=\xi''$ cannot achieve its
non-positive local minimum at a point in $(-\frac{\pi}{2},
\frac{\pi}{2})$. On the other hand, $\xi''(\pm\frac{\pi}{2})=2$,
by equation (\ref{xi-eq}), $\xi(\pm \frac{\pi}{2})=0$ and
$\xi'(\pm \frac{\pi}{2})=\pm \frac{2\pi}{3}$. Therefore
$\xi''(t)>0$ on $[-\frac{\pi}{2},\frac{\pi}{2}]$ and $\xi'$ is
increasing. Since $\xi'(t)=0$, we have $\xi'(t)< 0$ on
$(-\frac{\pi}{2},0)$ and $\xi'(t)>0$ on $(0,\frac{\pi}{2})$.
%satisfies the equation
%\[
%\frac12 h'' t\cos t + (\cos t - 2t \sin t) h' -2(\sin t + t\cos t)
%h = -4\sin t.
%\]
Similarly, from the equation
\begin{eqnarray}                                           \label{q''-eq}
&\frac{\cos^2 t}{2(1+\cos^2 t)}(q'')'' -\frac{\cos t\sin t
(3+2\cos^2 t)}{(1+\cos^2 t)^2}(q'')' -\frac{2(5\cos^2 t+\cos^4 t)}{(1+\cos^2 t)^2}(q'') \nonumber\\
&=-\frac{8\cos t\sin t}{(1+\cos^2 t)^2}
\end{eqnarray}
we get the results in the last line of the lemma.

Set $h(t)=\xi''(t)t-\xi'(t)$. Then $h(0)=0$  and $h'(t)=
\xi'''(t)t>0$ in $(0,\frac{\pi}{2})$. Therefore
$(\frac{\xi'(t)}{t})'=\frac{h(t)}{t^2}>0$ in $(0,\frac{\pi}{2})$.
Note that $\frac{\xi'(-t)}{-t}= \frac{\xi'(t)}{t}$,
$\frac{\xi'(t)}{t}|_{t=0}=\xi''(0)=2(3-\frac{\pi^2}{4})$ and
$\frac{\xi'(t)}{t}|_{t=\pi/2}=\frac43$. This completes the proof
of the lemma.
\end{proof}

\begin{lemma}                                       \label{eta-lemma}
Let
\begin{equation}                                    \label{eta-def}
 \eta(t)=\frac{\frac{4}{\pi}t+\frac{4}{\pi}\cos t\sin t-2\sin t}{\cos^2t}
 \qquad
\textrm{on}\quad [-\frac{\pi}{2},\frac{\pi}{2}\,].
\end{equation}
Then the function $\eta$ satisfies the following
\begin{eqnarray}
& &{}\frac{1}{2}\eta''\cos ^2t-\eta'\cos t\sin t-\eta
     =-\sin t\qquad \textrm{in \  }(-\frac{\pi}{2},\frac{\pi}{2}\,),              \label{eta-eq}\\
& &{}\eta'\cos t -2\eta \sin t =\frac{8}{\pi}\cos t -2  \qquad \textrm{in \  }(-\frac{\pi}{2},\frac{\pi}{2}\,),                         \label{eta-eq2}\\
& &{}-1=\eta(-\frac{\pi}{2})\leq\eta(t)\leq\eta(\frac{\pi}{2})=1 \qquad \textrm{on \ }[-\frac{\pi}{2},\frac{\pi}{2}\,],\nonumber\\
& &{}0<2(\frac{4}{\pi}-1)=\eta'(0)\leq \eta'(t)\leq
\eta'(\pm \frac{\pi}{2}) =\frac{8}{3\pi} \qquad \textrm{on \ }[-\frac{\pi}{2},\frac{\pi}{2}\,],\nonumber\\
& &{}-1/2=\eta''(-\frac{\pi}{2})\leq \eta''(t)\leq
\eta''(\frac{\pi}{2}) =1/2\qquad \textrm{on \ }[-\frac{\pi}{2},\frac{\pi}{2}\,],\nonumber\\
& &{}\eta'''(t)>0 \textrm{ \ on \ }[-\frac{\pi}{2},
\frac{\pi}{2}]\quad \textrm{and}\quad \eta'''(\pm
\frac{\pi}{2})=\frac{32}{15\pi}.\nonumber
\end{eqnarray}
\end{lemma}
\begin{proof}\quad Let $p(t)= \eta'(t)$, i.e.,
\begin{equation}                                                \label{p-def}
p(t) = \eta'(t) = \frac{2(\frac{4}{\pi}\cos t+\frac{4}{\pi}t\sin
t-\sin^2t-1)} {\cos^3t}.
\end{equation}
%Thus
%\[
%p'(t)=\eta''(t)=\frac{2[\frac{4}{\pi}t+\frac{8}{\pi}t\sin^2t+\frac{12}{\pi}\cos
%t \sin t-5\sin t-\sin^3 t]}{\cos^4 t}.
%\]
Equation (\ref{eta-eq}), $\eta(\pm \frac{\pi}{2})=\pm 1$,
$\eta'(0)=2(\frac{4}{\pi}-1)$ and $\eta'(\pm \frac{\pi}{2})
=\frac{8}{3\pi}$ can be verified directly.  We get $\eta''(\pm
\frac{\pi}{2}) =\pm 1/2$ from the above values and equation
(\ref{eta-eq}). By (\ref{eta-eq}), $q=\eta'$, $q'=\eta''$ and
$p''=\eta'''$ satisfy the following equations in
$(-\frac{\pi}{2},\frac{\pi}{2})$
\begin{equation}                                        \label{p-eq}
\frac12 p''\cos t -2p'\sin t -2p\cos t = -1,
\end{equation}
%By (\ref{p-eq}),
%\[
%\frac12p''-2p'\tan t -2p =-\sec t,
%\]
\[
\frac{\cos^2 t}{2(1+\cos^2 t)}p'''-\frac{2\cos t\sin t}{1+\cos^2
t}p''-2p'=-\frac{\sin t}{1+\cos^2 t},
\]
and
\begin{eqnarray}                                           \label{p''-eq}
&\frac{\cos^2 t}{2(1+\cos^2 t)}(p'')'' -\frac{\cos t\sin t
(3+2\cos^2 t)}{(1+\cos^2 t)^2}(p'')' -\frac{2(5\cos^2 t+\cos^4 t)}{(1+\cos^2 t)^2}(p'') \nonumber\\
&=-\frac{\cos t(2+\sin t)}{(1+\cos^2 t)^2}.
\end{eqnarray}
The coefficient of $(p'')$ in (\ref{p''-eq}) is obviously negative
in $(-\frac{\pi}{2}, \frac{\pi}{2})$ and the righthand side of
(\ref{p''-eq}) is also negative. So $p''$ cannot achieve its
non-positive local minimum at a point in $(-\frac{\pi}{2},
\frac{\pi}{2})$. On the other hand,
$p''(\frac{\pi}{2})=\frac{32}{15\pi}>0$ (see the proof below),
$p''(t)
> 0$ on $[-\frac{\pi}{2}, \frac{\pi}{2}]$. Therefore $p'$ is increasing
and $-1/2=p'(-\frac{\pi}{2})\leq p'(t)\leq p'(\frac{\pi}{2})
=1/2$. Note that $p'(0)= 0$ ($p'$ is an odd function). So $p'(t)>
0$ on $(0,\frac{\pi}{2})$ and $p$ is increasing on
$[0,\frac{\pi}{2}\,]$. Therefore $2(4/\pi-1)=p(0)\leq
p(t)=\eta'(t)\leq p(\frac{\pi}{2}) =\frac{8}{3\pi}$ on
$[0,\frac{\pi}{2}]$, and on $[-\frac{\pi}{2}, \frac{\pi}{2}]$
since $p$ is an even function.
%show that $p''(\frac{\pi}{2})=\frac{32}{15\pi}$
We now show that $p(\frac{\pi}{2})=\frac{8}{3\pi}$,
$p'(\frac{\pi}{2})=1/2$ and $p''(\frac{\pi}{2})=\frac{32}{15\pi}$.
The first is from a direct computation by using (\ref{p-def}). By
(\ref{eta-eq}),
\[
\frac12p'(\frac{\pi}{2})=\frac12\eta''(\frac{\pi}{2}) =
\lim_{t\rightarrow \frac{\pi}{2}^-}\frac{\eta'(t)\cos t\sin t +
\eta(t)-\sin t}{\cos^2 t}=-\frac12[\eta''(\frac{\pi}{2})-1].
\]
So $p'(\frac{\pi}{2})=1/2$. Similarly, by (\ref{p-eq}),
\[
\frac12p''(\frac{\pi}{2}) =\lim_{t\rightarrow
\frac{\pi}{2}-}\frac{2p'(t)\sin t-1}{\cos t}  +
2p(\frac{\pi}{2})=-2p''(\frac{\pi}{2}) + \frac{16}{3\pi}
\]
Thus $p''(\frac{\pi}{2})=\frac{32}{15\pi}$.
%end---show that $p''(\frac{\pi}{2})=\frac{32}{15\pi}$
\end{proof}

\begin{lemma}                                                   \label{r-lemma}
The function \ $r(t)=\xi'(t)/\eta'(t)$ is an increasing function
on $[-\frac{\pi}{2}, \frac{\pi}{2}]$, i.e., $r'(t)>0$, and $|r(t)|
\leq \frac{\pi^2}{4}$ holds on $[-\frac{\pi}{2}, \frac{\pi}{2}]$.
\end{lemma}
\begin{proof}\quad Let $p(t) =\eta'(t)$ as in (\ref{p-def})
and $q(t)=\xi'(t)$. Then $r(t)=q(t)/p(t)$. It is easy to  verify
that $ r(\pm \frac{\pi}{2})=\pm \frac{\pi^2}{4}$. By (\ref{p-eq})
and (\ref{q-eq}),
\[ (1/2)p(t)r''\cos t +(p'(t)\cos t-2p(t)\sin t)r'-r=-4\sin t.
\]
Differentiating the last equation, we get
\begin{eqnarray}
 &[\frac12p(t)\cos t] (r')''+[\frac32p'(t)\cos t-\frac52 p(t)\sin t]
(r')'\nonumber\\
 &+[p''(t)\cos t -3p'(t)\sin t - 2p(t)\cos t -1](r')=-4\cos t.\nonumber
\end{eqnarray}
Using (\ref{p-eq}), the above equation becomes
\begin{eqnarray}                                       \label{r'-eq}
 &[\frac12p(t)\cos t] (r')''+[\frac32p'(t)\cos t-\frac52 p(t)\sin t]
(r')'\nonumber\\
 &+[p'(t)\sin t + 2p(t)\cos t -3](r')=-4\cos t.
\end{eqnarray}
The coefficient of $(r')$ in (\ref{r'-eq}) is negative, for
$p'(t)\sin t +2p\cos t -3< \frac12+ \frac{16}{3\pi} -3 < 0$. This
fact and the negativity of the righthand side of (\ref{r'-eq}) in
$(-\frac{\pi}{2},\frac{\pi}{2})$ imply that $r'$ cannot achieve
its non-positive minimum on $[-\frac{\pi}{2},\frac{\pi}{2}]$ at a
point in $(-\frac{\pi}{2}, \frac{\pi}{2})$.  Now
\begin{eqnarray}\nonumber
  & &\lim_{t\rightarrow \frac{\pi}{2}^-} r'(t) \nonumber\\
  & =&\lim_{t\rightarrow \frac{\pi}{2}^-}s(t)\cos^2 t
  /(\frac{4}{\pi}\cos t +\frac{4}{\pi}t\sin t- \sin^2 t
-1)^2\nonumber\\
 &=&\lim_{t\rightarrow \frac{\pi}{2}^-}[s(t)/ \cos^4]
  /[(\frac{4}{\pi}\cos t +\frac{4}{\pi}t\sin t- \sin^2 t
-1)/ \cos^3 t]^2\nonumber\\
 &=&\lim_{t\rightarrow \frac{\pi}{2}^-}[ s(t)/\cos^4 t]
  /[\frac12\eta'(t)]^2\nonumber\\
  &=&(\frac{4}{3\pi}-\frac{\pi}{12})/(\frac{4}{3\pi})^2\nonumber\\
  &>&0,\nonumber
\end{eqnarray}
where
\begin{eqnarray}\nonumber
s(t)&=&-\frac{4}{\pi}t^2-t^2 \cos t +\frac{12}{\pi}\cos^2 t
  +\frac{8}{\pi}t\sin t\cos t \nonumber\\
  & &-\cos t\sin^2t
 +(\frac{\pi^2}{4}-3)\cos t - \pi + 4t\sin t.\nonumber
\end{eqnarray}
Therefore $r'(t)> 0$ and $r$ is an increasing function on
$[-\frac{\pi}{2},\frac{\pi}{2}]$.
\end{proof}

\begin{proof}[Proof of the estimate (\ref{delta-def})]
We estimate the maximum of the function
\begin{equation}                \label{p-of-x-def}
P(x)=|\nabla v|^2+Av^2,
\end{equation}
where $v$ is the function in (\ref{v-def}), and where $A\geq 0$ is
a constant.

Let $A=0$ in (\ref{p-of-x-def}). Function P must achieve its
maximum at some point $x_0\in \bar{\Omega}$. Suppose that
$x_0\in\partial\Omega$. Choose an orthornormal frame $\{ e_1,
\dots, e_n\}$ about $x_0$ such that $e_n$ is  a outward normal to
$\partial \Omega$. By (\ref{basic2.001}), $v_n=\partial v/\partial
N=0$. Thus at $x_0$
\begin{eqnarray}
v_{in}&=&e_ie_nv-(\nabla _{e_i}e_n)v\nonumber\\
      &=&-(\nabla _{e_i}e_n)v\nonumber\\
      &=&-\sum_{j=1}^{n-1} h_{ij}v_j\nonumber
\end{eqnarray}
and
\begin{eqnarray}
P_n&=&2\sum_{j=1}^nv_jv_{jn}+2Avv_n=\sum_{j=1}^{n-1}v_jv_{jn}\nonumber\\
      &=&-2\sum_{i, j=1}^{n-1} h_{ij}v_iv_j\nonumber\\
      &\leq& 0 \qquad \textrm{by the convexity of }\partial \Omega.\nonumber
\end{eqnarray}
On the other hand, that $P$ attains the maximum at the boundary
point $x_0$ implies that
\[
P_n\geq 0.
\]
Thus at $x_0$, $-2\sum_{i, j=1}^{n-1} h_{ij}v_iv_j=P_n=0$. By the
strict convexity, $v_1=\dots=v_{n-1}=0$ and $\nabla v=0$ at $x_0$.
Therefore $v$ is a constant. This is impossible, so $x_0\in
\Omega$. $\nabla v (x_0)\not=0$ (otherwise $v$ is a constant). At
$x_0$,
\[
\nabla P(x_0)=0 \qquad \textrm{and}\qquad \Delta P(x_0)\leq 0.
\]
Take a local frame so that
\[
v_1(x_0)=\nabla v(x_0)\qquad \textrm{and}\qquad v_i(x_0)=0,\quad
i\geq2.
\]
Thus at $t_0$ we have
\[
0=\frac12\nabla P_i=v_jv_{ji}+Avv_i,
\]
\begin{equation}                \label{d1}
v_{11}=-Av \qquad \textrm{and}\qquad v_{1i}=0\quad i\geq 2,
\end{equation}
and
\begin{eqnarray}
&0 &\geq \frac12\Delta P(x_0) =v_{ji}v_{ji}+v_{j}v_{jii}+Av_{i}v_{i}+Avv_{ii}\nonumber\\
&{} &=v_{ji}^2+v_1(\Delta v)_1 +R_{ji}v_{j}v_{i}+A|\nabla v|^2 +Av\Delta v\nonumber\\
&{} &\geq v_{11}^2+v_1(\Delta v)_1 +A|\nabla v|^2 +Av\Delta v\nonumber\\
&{} &=(-Av)^2-\lambda |\nabla v|^2-2v_{1}(\nabla v\nabla \ln f )_{1}+A|\nabla v|^2\nonumber\\
&{} &{ } -\lambda Av(v+a)-2Av\nabla v\nabla \ln f\nonumber\\
&{} &=-(\lambda-A)|\nabla v|^2-Av^2(\lambda-A)-a\lambda A v\nonumber\\
&{} &{} -2v_1^2(\ln f)_{11}-2v_1(\ln f)_1(v_{11}+A v),\nonumber
\end{eqnarray}
where we have used (\ref{d1}) and (\ref{delta-def}). Therefore at
$x_0$,
\begin{equation}
0\geq -(\lambda -2\alpha -A)|\nabla v|^2-A(\lambda-A)v^2-a\lambda
A v.                     \label{d2}
\end{equation}
Using the fact that $A=0$ in the above inequality, we get the
(\ref{delta-bound}).
\end{proof}

\begin{proof}[Proof of the estimate (\ref{basic5})]
We first prove the following
\begin{equation}                        \label{basic3}
\frac{\left |\nabla v\right |^2}{b^2-v^2} \leq\lambda(1+a),
\end{equation}
where $b>1$ is an arbitrary constant.

Let $A=\lambda (1+a)+\epsilon$ in (\ref{p-of-x-def}) for small
$\epsilon>0$. $P$ achieves its maximum at some
$x_0\in\bar{\Omega}$. If $\nabla v(x_0)\not=0$ and $x_0\in\Omega$,
then (\ref{d2}) holds at $x_0$ with $A=\lambda (1+a)+\epsilon$.
Thus
\[
|\nabla v(x_0)|^2+ \lambda(1+a)v(x_0)^2\leq \frac{a\lambda
v}{a\lambda +\epsilon}[\lambda(1+a) +\epsilon]\leq [\lambda(1+a)
+\epsilon].
\]
This estimate holds if $x_0\in \bar{\Omega}$ with $\nabla
v(x_0)=0$. If $x_0\in\partial\Omega$, then the convexity of
$\Omega$ and previous argument in the proof of (\ref{delta-bound})
imply that the above estimate holds. So we have the estimate
(\ref{basic3}). By the definition of $Z$, we have (\ref{basic5}).
\end{proof}

\begin{proof}[Proof of Lemma \ref{barrier-lemma}]\quad Define
\[ J(x)=\left\{ \frac{\left |\nabla v\right |^2}{b^2-v^2}
-\lambda z \right\}\cos^2t,
\]
where $t=\sin^{-1}(v(x)/b)$. Then
\[ J(x)\leq 0\quad\textrm{for } x\in \bar{\Omega}
\qquad \textrm{and} \qquad J(x_0)=0.
\]
If $\nabla v(x_0)=0$ then
\[ 0=J(x_0)=-\lambda z\cos^2 t.
\]
This contradicts the condition 3 in the theorem. Therefore
\[ \nabla v(x_0)\not=0.
\]
If $x_0\in \partial \Omega$, then by an argument in the proof of
(\ref{delta-bound}), the convexity of $\Omega$ and that $J(x_0)$
is the maximum would imply that $\nabla v(x_0)=0$. Thus $x_0\in
\Omega =\bar{\Omega}\backslash\partial \Omega$. The Maximum
Principle implies that
\begin{equation}                                                \label{es1}
\nabla J(x_0)=0\qquad \textrm{and}\qquad \Delta J(x_0)\leq 0.
\end{equation}
$J(x)$ can be rewritten as
\[  J(x)=\frac{1}{b^2}|\nabla v|^2-\lambda z\cos^2t.
\]
Thus (\ref{es1}) is equivalent to
\begin{equation}                                              \label{es2}
\frac{2}{b^2}\sum_{i}v_iv_{ij}\Big|_{x_0}=\lambda\cos t[z' \cos t
-2z\sin t]t_j\Big|_{x_0}
\end{equation}
and
\begin{eqnarray}                                              \label{es3}
0&\geq&\frac{2}{b^2}\sum_{i,j}v_{ij}^2+\frac{2}{b^2}\sum_{i,j}v_iv_{ijj}
 -\lambda (z''|\nabla t|^2+z'\Delta t)\cos^2t \\
 & &+4\lambda z'\cos t\sin t |\nabla t|^2 -
\lambda z\Delta\cos^2t\Big|_{x_0}.\nonumber
\end{eqnarray}
Choose a normal coordinate around $x_0$ such that $v_1(x_0)\not=0$
and $v_i(x_0)=0$ for $i\geq 2$. Then (\ref{es2}) implies
\begin{equation}                                             \label{es4}
v_{11}\Big|_{x_0}=\frac{\lambda b}{2}(z'\cos t-2z\sin t)
\Big|_{x_0}\quad\textrm{and}\quad v_{1i} \Big|_{x_0}=0\
\textrm{for } i\geq2.
\end{equation}
Now we have
\begin{eqnarray}
|\nabla v|^2
\Big|_{x_0}&=&\lambda b^2z\cos^2t\Big|_{x_0},\nonumber\\
 |\nabla t|^2
\Big|_{x_0}&=&\frac{|\nabla v|^2}{b^2-v^2}=\lambda z
\Big|_{x_0},\nonumber\\
\frac{\Delta v}{b}\Big|_{x_0} &=&\Delta \sin t =\cos t\Delta
t-\sin t |\nabla t|^2
\Big|_{x_0},\nonumber\\
\Delta t\Big|_{x_0}&=&\frac{1}{\cos t}(\sin t|\nabla
t|^2+\frac{\Delta v}{b})
\nonumber\\
 &=&\frac{1}{\cos t} [ \lambda z\sin t-\frac{\lambda}{b}(v+a)
-\frac{2}{b}v_1(\ln f)_1] \Big|_{x_0}, \quad\textrm{and}
\nonumber\\
\Delta\cos^2t\Big|_{x_0}&=&\Delta \left(1-\frac{v^2}{b^2}\right)
 =-\frac{2}{b^2}|\nabla v|^2-\frac{2}{b^2}v\Delta v
\nonumber\\
  &=&-2\lambda z\cos^2t+\frac{2}{b^2}\lambda v(v+a)
+\frac{4}{b^2}vv_1(\ln f)_1\Big|_{x_0}. \nonumber
\end{eqnarray}
Therefore,
\begin{eqnarray}
& {}&
\frac{2}{b^2}\sum_{i,j}v_{ij}^2\Big|_{x_0}\geq\frac{2}{b^2}v_{11}^2
\nonumber\\
& {}& =\frac{\lambda ^2}{2}(z')^2\cos^2t-2\lambda ^2zz'\cos t\sin
t
      +2\lambda ^2z^2\sin^2 t\Big|_{x_0}\nonumber,
\end{eqnarray}
\begin{eqnarray}
\frac{2}{b^2}\sum_{i,j}v_iv_{ijj}\Big|_{x_0}
&=&\frac{2}{b^2}\left(\nabla v\,\nabla
       (\Delta v)+R(\nabla v,\nabla v)\right)\
\geq \frac{2}{b^2}\nabla v\,\nabla (\Delta v)\nonumber\\
 &=&-2\lambda^2z\cos^2t-\frac{4}{b^2}v_1v_{11}(\ln f)_1
-\frac{4}{b^2}v_1^2(\ln f)_{11}\Big|_{x_0},\nonumber
\end{eqnarray}
\begin{eqnarray}
&{}&  -\lambda (z''|\nabla t|^2+ z'\Delta t)\cos^2t\Big|_{x_0}\nonumber\\
&{}&=-\lambda^2 zz''\cos^2t-
\lambda^2zz'\cos t\sin t\nonumber\\
&{}&{ }+\frac{1}{b}\lambda^2z'(v +a)\cos t +\frac{2}{b}\lambda
z'v_1(\ln f)_1\cos t\Big|_{x_0},\nonumber
\end{eqnarray}
and
\begin{eqnarray}
&{}&4\lambda z'\cos t\sin t|\nabla t|^2-\lambda z\Delta
\cos^2t\Big|_{x_0}
\nonumber\\
&{}&=4\lambda^2zz'\cos t\sin t+2\lambda^2z^2\cos^2t\nonumber\\
&{}&{ }-\frac{2}{b}\lambda^2z\sin t\,(v+a) -\frac{4}{b}\lambda
z\sin t\,v_1(\ln f)_1\Big|_{x_0}.\nonumber
\end{eqnarray}
Putting these results into (\ref{es3}) we get
\begin{eqnarray}                                                   \label{es5}
0&\geq&-\lambda^2zz''\cos^2t+ \frac{\lambda^2}{2}(z')^2\cos^2t
\nonumber\\
 & & {}  +\lambda^2z'\cos t\left(z\sin t+c +\sin t\right)
  \\
 & & {}+2\lambda^2z^2-2\lambda^2z\nonumber\\
 & & {}-2\lambda^2cz\sin t -4\lambda z\cos^2t(\ln f)_{11}
\nonumber\\
 & & {}- \frac{4}{b^2}\left[v_{11}-\frac{\lambda b}{2}(z'\cos t
-2z\sin t)  \right]v_1(\ln f)_1
 \Big|_{x_0}.\nonumber
\end{eqnarray}
The last term in (\ref{es5}) is 0 due to (\ref{es4}) . Now
\begin{equation}                                                    \label{es6}
z(t_0)>0,
\end{equation}
by the condition 3 in the theorem, and
\begin{equation}                                                    \label{es9}
-\frac{(\ln f)_{11}}{\lambda}\geq \delta,
\end{equation}
by the definition of $\delta$. Dividing two sides of (\ref{es5})
by $2\lambda^2z\Big|_{x_0}$ and taking (\ref{es9}) into account,
we have
\begin{eqnarray}
0&\geq&-\frac12z''(t_0)\cos^2t_0 +\frac12z'(t_0)\cos t_0\left(\sin
t_0
            +\frac{c +\sin t_0}{z(t_0)}\right) +z(t_0) \nonumber\\
 & & {}  -1-c\sin t_0 +2\delta \cos^2t_0\nonumber\\
 & & {}+\frac{1}{4z(t_0)}(z'(t_0))^2\cos^2t_0.\nonumber
\end{eqnarray}
Therefore,
\begin{eqnarray}
0&\geq&-\frac12z''(t_0)\cos^2t_0 + z'(t_0)\cos t_0\sin t_0+z(t_0)
            -1-c\sin t_0 +2\delta \cos^2t_0\nonumber\\
 & & {}+\frac{z'(t_0)}{4z(t_0)}\cos t_0[z'(t_0)\cos t_0 -2z(t_0)\sin t_0 +
            2\sin t_0 + 2c].\nonumber
\end{eqnarray}
\end{proof}

\begin{proof}[Proof of Corollary  \ref{corollary1}] \quad By Condition 2
in the theorem, (\ref{basic5}), $|\sin t_0|=|v(t_0)/b|\leq1/b$ and
$1-c \leq z(t_0)\leq 1+a$. Thus for $t_0\geq 0$,
\[
-z(t_0)\sin t_0+ \sin t_0 + c  \geq -\sin t_0-a\sin t_0+\sin
t_0+c\geq a(\frac{1}{b}-\sin t_0)\geq 0,
\]
and for $t_0<0$,
\[
-z(t_0)\sin t_0+ \sin t_0 + c  \geq -\sin t_0 + c\sin t_0 +\sin
t_0 +c\geq c(1+\sin t_0)\geq 0.
\]
In any case the last term in the (\ref{barrier-eq}) is
non-negative.
\end{proof}

\begin{proof}[Proof of Corollary  \ref{corollary2}]\quad The last term
in the (\ref{barrier-eq}) is nonnegative.
\end{proof}

The following Lemma is due to the author \cite{ling1}. We enclose
it here for the completeness.
\begin{lemma}                   \label{convex-bound-lemma2}
Let $f$ be the first eigenfunction of (\ref{basic1}) with $f>0$ in
$\Omega$. Then there exists an $\epsilon>0$ such that the function
$-\ln f$ is strictly convex in the $\epsilon$-neighborhood of
$\partial\Omega$.
\end{lemma}
\begin{proof}\quad Choose a normal coordinate about $x_0$ such
that $\partial/\partial x_1$ is the outward unit normal vector
field of $\partial \Omega$ near $x_0$. Take a point $\bar x$ with
small distance $d$ to $\partial \Omega$ and $d=dist(\bar x,
x_0)=dist(\bar x, \partial \Omega)$. Then by the Strong Maximum
Principle $f_1|_{x_0}<0$ , and $f_i|_{x_0}=0$, for $i\geq 2$.
Therefore
\[
f_1|_{\bar x} \sim c_1 d,\quad c_1=-f(x_0)>0
\]
and
\[
f_i|_{\bar x}\sim O(d)\quad\textrm{for}\quad i\geq 2.
\]
Here "$\sim A/d^\alpha$" means "$=(A+o(1))/d^\alpha$" with
$o(1)\rightarrow 0$ as $d\rightarrow 0$. Let $w=\ln f$. Then
\[
w_{ij}\Big |_{\bar x} =(f_{ij}/f-f_if_j/f^2-\Gamma_{ij}^k
f_k/f]\Big |_{\bar x}.
\]
Thus we obtain
\[
w_{11}\Big |_{\bar x}\sim O(1/d)-f_1^2/d^2
\]
and
\[w_{1i}\Big |_{\bar x}\sim O(1/d),\quad i\geq 2.
\]
For $i, j\geq 2$,
\[
f_{ij}\Big |_{\bar x}=\nabla_{\frac{\partial}{\partial x_i}}
\nabla_{\frac{\partial}{\partial x_j}}f\Bigg |_{\bar x}
-\nabla_{\frac{\partial}{\partial x_i}} \frac{\partial}{\partial
x_j}f\Bigg |_{\bar x} \sim O(d) -\nabla_{\frac{\partial}{\partial
x_i}} \frac{\partial}{\partial x_j}f\Bigg |_{x_0}.
\]
Now
\[\nabla_{\frac{\partial}{\partial x_i}}
\frac{\partial}{\partial x_j}f\Bigg |_{x_0}
=(\nabla_{\frac{\partial}{\partial x_i}} \frac{\partial}{\partial
x_j},\frac{\partial}{\partial x_k}) \frac{\partial}{\partial
x_k}f\Bigg |_{x_0} =(\nabla_{\frac{\partial}{\partial x_i}}
\frac{\partial}{\partial x_j},\frac{\partial}{\partial x_1})
\frac{\partial}{\partial x_1}f\Bigg |_{x_0},
\]
\[
(\frac{\partial}{\partial x_1},\frac{\partial}{\partial x_j})=0,
\quad j\geq 2,
\]
and
\[
0=\frac{\partial}{\partial x_i} (\frac{\partial}{\partial
x_1},\frac{\partial}{\partial x_j})
=(\nabla_{\frac{\partial}{\partial x_i}} \frac{\partial}{\partial
x_1},\frac{\partial}{\partial x_j}) +(\frac{\partial}{\partial
x_1},\nabla_{\frac{\partial}{\partial x_i}}
\frac{\partial}{\partial x_j}) =h_{ij}
+(\nabla_{\frac{\partial}{\partial x_i}} \frac{\partial}{\partial
x_j},\frac{\partial}{\partial x_1}),
\]
where $(h_{ij})_{n-1,n-1}$ is the second fundamental form of
$\partial\Omega$ to $\partial/\partial x_1$. Therefore
\[
f_{ij}\Big |_{\bar x}\sim O(d)+h_{ij}f_1\Big |_{x_0}
\quad\textrm{and}\quad w_{ij}\Big |_{\bar x}\sim h_{ij}f_1\Big
|_{x_0},\quad i,j\geq 2,
\]
and
\[
(-w_{ij})_{n,n}\Big |_{\bar x} \sim \left(
\begin{array}{cc}
f_1^2/d^2  & O(1/d) \\
O(1/d) &  O(1/d)(-f_1h_{ij})_{n-1,n-1}
\end{array}
\right)\Bigg |_{x_0}
\]
Since $(h_{ij})$ is positive definite, so is $(-w_{ij})\Big
|_{\bar x}$ for $\bar x$ near the boundary $\partial \Omega$.
\end{proof}

%%%%%%%%%%%%End of Real Paper%%%%%%%%%%%
%
%
%%%%%%%%%%%%%%%%%References%%%%%%%%%%%%%%%%%%%
%\nocite{*}
%\bibliographystyle{plain}
%\bibliography{GapBound} %GapBound.bib, for using BibTeX
%\input{gapBound-sub.bbl}

Department of mathematics, Utah Valley State College, Orem, Utah
84058

\textit {E-mail address}: \texttt{lingju@uvsc.edu}
\end{document}